\documentclass{article}

\usepackage[cp1251]{inputenc}
\usepackage[english,russian]{babel}
\usepackage[tbtags]{amsmath}
\usepackage{amsfonts,amssymb}

\usepackage{graphicx}

\newtheorem{theorem}{Теорема}

\newtheorem{lemma}[theorem]{Лемма}

\newtheorem{cor}[theorem]{Следствие}

\begin{document}


\author{
В.\,Л.~Гейнц, А.\,А. ~Шкаликов
}

\title{
\begin{center}
Оценка отношения   двух целых функций, \\
нули которых совпадают в круге.
\end{center}
}


\maketitle

\begin{abstract}
 В  заметке изучаются целые функции конечного порядка роста, которые на некотором луче комплексной
плоскости допускают представление $\psi(z) = 1+ O(|z|^{-\mu}),  \mu >0$.
Получен следующий результат: если нули двух функций $\psi_1, \psi_2$  такого класса совпадают в круге
радиуса  $R$ с центром в нуле,
то при любых произвольно малых $\delta\in (0,1), \ \varepsilon >0$
в круге радиуса $R^{1-\delta}$  отношение этих функций допускает оценку
$|\psi_1(z)/\psi_2(z) -1| \leqslant \varepsilon R^{-\mu(1-\delta)}$,  если $R\geqslant R_0(\varepsilon, \delta)$
Полученны результаты важны для анализа устойчивости
в задаче о восстановлении потенциала уравнения Шредингера на полуоси по резонансам оператора.

Библиография: 3 названия.
\end{abstract}

\footnotetext{Работа выполнена при финансовой поддержке гранта РФФИ 15-51-12389.}

\bigskip

{\bf 1. Введение.} В спектральной теории оператора Шредингера с быстроубывающем потенциалом важную роль играют целые функции
вида
\begin{equation}\label{0}
\psi(z) = 1 +\int_0^\infty K(t) e^{izt}\, dt,
\end{equation}
где функция $K(t)$ определяется ядром оператора преобразования, которое в свою
очередь определяется потенциалом  $q$ оператора Шредингера
$$
A y =-y'' +q(x)y, \qquad x\in [0,\infty ).
$$
Известно \cite[Гл. 2]{Ma}, что если потенциал $q$  подчинен условию
$$
\int_0^\infty |q(x)|exp (x^\gamma)\, dx \leqslant C_0 <\infty, \quad \text{ при некотором }\ \, \gamma >1,
$$
то  функция $K$  подчинена оценке
\begin{equation*}
|K(t)|\leqslant C\exp\left(-\left(\frac{t}{2}\right)^{\gamma}\right),
\end{equation*}
с постоянной $C$, зависящей от $C_0$.
В этом случае нетрудно показать, что
функция $\psi$ является целой порядка роста, не превосходящего
\begin{equation*}
\rho(\gamma) = \frac{\gamma}{\gamma - 1}
\end{equation*}
и конечного типа. Кроме того, при $z\in \mathbb{C^+}$ выполняется оценка \cite[Лемма 3.1.3]{Ma}
\begin{equation}\label{im z}
|\psi(z) - 1|\le \frac{C_1}{\text{Im} z}, \qquad \text{Im} z >0,
\end{equation}
с постоянной $C_1$,   зависящей от  $C$  и $\gamma$ .
В случае, если потенциал $q$ финитный, функция $K$
также оказывается финитной, а тогда
функция \eqref{0} является целой функцией экспоненциального типа, т.е. порядка 1 и конечного типа.

Нули функции $\psi$ называются резонансами оператора Шредингера, которые
 являются физически наблюдаемыми величинами.
Поэтому  важное значение имеет задача о восстановлении  потенциала $q$ по
резонансам оператора Шредингера.
 В решении этой задачи ключевую роль играют оценки отношения двух целых функций вида \eqref{0} при условии,
 что нули этих функций совпадают в круге большого радиуса $R$.  Важность и эффективность
 таких оценок продемонстрирована в работе Марлетты, Штеренберга и Вайкарта \cite{MaSW}.

 Основная цель нашей работы --- получить обобщение и усиление  одного из ключевых
 аналитических результатов работы \cite{MaSW}. Применение полученных здесь результатов авторы
 планируют провести в другой работе.

\bigskip
{\bf 2. Основной результат.}  Обозначим через $B(a,r)$ открытый круг радиуса $r$ с центром в точке $a\in \mathbb C$, а через
$\mathcal{M} = \mathcal{M}(C_0,C_1, \rho, \sigma, \mu, r_0)$
множество, состоящее из целых функций $\psi$, удовлетворяющих следующим условиям:

1) функции $\psi \in \mathcal M$ подчинены оценке
\begin{equation}\label{1}
|\psi(z)| \leqslant C_0 e^{\sigma r^\rho}   \quad\text{для всех} \  {|z|\geqslant r_0}\geqslant 1;
\end{equation}

2) На некотором луче $\arg z = \varphi \in [0, 2\pi)$  справедлива оценка
\begin{equation}\label{2}
|\psi(z) -1| \leqslant \frac{C_1}{|z|^\mu}, \quad  \mu>0, \qquad \text{для всех}\ \,
|z|\geqslant r_0.
\end{equation}
Значение аргумента $\varphi$ далее не играет роли и мы для упрощения записи полагаем $\varphi =0$.
Все постоянные, участвующие в определении класса $\mathcal M$,  предполагаются положительными.
Далее через $[\rho]$  обозначается целая часть числа $\rho$.

Для функций из $\mathcal M $  мы получим следующий результат.


\begin{theorem}.
Пусть функции $\psi_1, \psi_2 \in \mathcal M$ и нули этих функций в круге $B(0,R)$ совпадают с учетом кратностей.
Тогда при любых  произвольно малых $\delta\in (0,1)$ и $\varepsilon >0$  для отношения функций $\psi_1/\psi_2$
в круге $B(0,R^{1-\delta})$ справедлива оценка
\begin{equation}\label{result}
\left|\frac{\psi_2(z)}{\psi_1(z)} - 1\right| < \frac \varepsilon {R^{\, \mu(1-\delta)}},
\end{equation}
при условии, что $R\geqslant R_0$,  где $R_0$ достаточно большое число.  Постоянная $R_0$  зависит от
 $\delta, \varepsilon$ и всех параметров
 $C_0, C_1, \sigma, \rho, \mu, r_0$,  определяющих класс $\mathcal M$,  но не зависит от
 функций $\psi_1, \psi_2\ \in \mathcal M$.  Зависимость от параметров  можно явно указать.
 А именно, можно положить
  $$
 R_0 = \max \left(r_1, r_2, r_3, r_4, r_5, R'(\varepsilon, \delta) \right),
 $$
 где величина $ R'(\varepsilon, \delta)$ определена в \eqref{R_0}, а величины $r_j$  опредляются в  процессе
 доказательства в соотношениях \eqref{r 1}, \eqref{r 2}, \eqref{r 3}, \eqref{r 4},
  и  \eqref{r 5}. При этом в \eqref{r 2}  нужно положить $a=p+1$, а число $p\in \mathbb N$  определяется условием
  $$
  (\mu +\rho)/\delta \leqslant p+1 <  1 + (\mu +\rho)/\delta.
  $$
  \end{theorem}


{\bf Доказательство}. {\it Шаг 1.}
Сначала заметим, что достаточно доказать оценку \eqref{result}  с произвольной постоянной $C=C(\delta)$
в числителе вместо $\varepsilon$. Действительно, если такая оценка уже доказана, для произвольного $\delta\in (0,1)$
с постоянной $C=C(\delta)$ и $R_0 =R_0(\delta)$,  то достаточно взять число $\delta_1 \in (0, \delta)$
(например, $\delta_1 = \delta/2$) и заметить, что левая часть в \eqref{result} будет оцениваться величиной
\begin{equation*}
\frac {\tilde C_1}{R^{1-\delta_1}} = \frac {\tilde C_1}{R^{1-\delta}\, R^{\delta -\delta_1}} \leqslant
\frac\varepsilon{R^{1-\delta}},
\end{equation*}
если
\begin{equation}\label{R_0}
 R\geqslant R'_0(\varepsilon, \delta)= \max \left( R_0(\delta_1),
\left(\frac{\tilde C_1}\varepsilon\right)^{1/(\delta -\delta_1)}\right),
\quad \text{ где}\ \, \tilde C_1 = C(\delta_1).
\end{equation}
Итак, далее докажем оценку  \eqref{result}  при  $R\geqslant R_0$  с постоянной $C$  вместо $\varepsilon$.
При этом постоянные $C$ и $R_0$
 предъявим в явном виде в зависимости от $\delta$  и параметров класса $\mathcal M$. А именно,
 $$
 C= 20\, A_p\, C_1,  \quad R_0 =\max (r_1, r_2, r_3, r_4, r_5),
 $$
 где $C_1$ --- постоянная из оценки \eqref{2}, число $A_p$ определено \eqref{A p}, а числа $r_j$ явно указаны в ходе доказательства.

\medskip

  Наметим дальнейший план доказательства. Возможность получить нужную оценку основана на
 представлении \eqref{represent},  где бесконечные произведения $\Pi_1, \Pi_2$ учитывают только те нули
функций $\psi_1, \psi_2$,  которые   лежат вне круга $|z|\leqslant R$.  Мы покажем, что функция $\Phi =\Pi_2/\Pi_1 -1$
в круге $B\left(0, a R^{1-\delta}\right)$ допускает оценку $\leqslant CR^{-\alpha}$, где $\alpha =
\delta(p+1)-\rho$, $C$ --  постоянная и $a$ --- свободный параметр.
Выбор числа $p\in\mathbb N$ в нашей власти; возьмем наименьшее $p$,
 при котором $\alpha \geqslant \mu$ (в действительности, достаточно выполнения условия
   $\alpha \geqslant \mu(1-\delta)$,  но технически проще считать  $\alpha \geqslant \mu$). Из определения функций
    класса $\mathcal M$  следует, что оценка \eqref{result}  выполняется на положительном луче.
    Но тогда из \eqref{represent} следует, что функцию $h= e^{g_2-g_1} -1$  можно  оценить
     на отрезке $\Delta =[r, (p+1)r]$  при $r= a(p+1)^{-1}  R^{1-\delta}$. Далее мы покажем, что оценка
     функции $h$  на отрезке  $\Delta$ влечет ее оценку в круге $B(0, ar)$. Но тогда  из представления \eqref{represent2}
     получаем аналогичную оценку функции $\psi_2 /\psi_1 -1$  в том же круге $B(0, ar)$.
     Отметим, что за счет выбора свободного параметра $a$  можно отслеживать
     постоянные в оценках (и делать их более оптимальными).  Однако наиболее простая форма записи получается при
     $a=p+1$,  поэтому далее мы делаем именно этот выбор. Приступим к реализации намеченной цели.

\medskip

{\it Шаг 2.}
Обозначим через $n(a,r)$ число нулей  функции $\psi$  в круге радиуса $r$  с центром в $a$.
При $a=0$ будем писать $n(r)$. Оценим $n(r)$ для функций $\psi \in \mathcal M$. Из формулы Йенсена и оценки
\eqref{1} легко получается оценка $n(r) <  2\sigma (er)^\rho$  асимптотически при $r>R_0$. Однако $R_0$
зависит от $\psi$. Но при дополнительном условии \eqref{2} эту зависимость можно  устранить.

Положим
\begin{equation}\label{c}
c:= \max \left(r_0, (2C_1)^{1/\mu}\right),
\end{equation}
где $r_0, C_1$ ---  постоянные, участвующие в оценке \eqref{2}. Тогда
$$
|\psi(c)| >1 -1/2 =1/2,   \quad\text{поэтому} \ \, -\ln |\psi(c)| < \ln 2.
$$
Согласно формуле Йенсена с учетом оценки \eqref{1} получаем
\begin{multline*}
n(r)= n(0,r) \leqslant n(c, n+c) \leqslant
\int_{r+c}^{e(r+c)} \frac{n(c,t)}t\, dt \leqslant  \int_0^{e(r+c)} \frac{n(c,t)}t\, dt \\ =
\frac{1}{2\pi}\int_0^{2\pi}\ln |\psi(c+e(r+c)e^{i\varphi})|\, d\varphi - \ln |\psi(c)|\\ \leqslant
\ln C_0 +\sigma(e(r+c))^\rho + \ln 2 \leqslant \ln(2C_0) +\sigma(2e)^\rho r^\rho.
\end{multline*}
Следовательно,
\begin{equation}\label{n(r)}
n(r) \leqslant  2\sigma (2e)^\rho r^\rho,
\end{equation}
если только $r\geqslant c$ и $\ln(2C_0) \leqslant \sigma(2e)^\rho r^\rho$ или
\begin{equation}\label{r 1}
r\geqslant r_1, \quad\text{где} \ \ r_1 = \max \left( r_0,\ (2C_1)^{1/\mu},\
\frac 1{2e} \left(\frac{\ln(2C_0)}\sigma\right)^{1/\rho} \right).
\end{equation}

\medskip

{\it Шаг 3}.
При любом $p\geqslant [\rho]$ согласно   теореме Адамара \cite[ Гл. 4.2]{Le} функция $\psi \in \mathcal{M}$ допускает представление
\begin{equation*}
\psi(z) = z^{n}\exp(g(z))\prod\limits_{n=1}^{\infty}E_{p}\left(\frac{z}{z_n}\right),
\end{equation*}
где $p \ge [\rho]$,
\begin{equation*}
E_p(z) = (1-z)\exp\left(\sum\limits_{k=1}^{p} \frac{z^k}{k}\right),
\end{equation*}
 $g$--- многочлен степени $\leqslant p$ и $\{z_n\}_{n=1}^{\infty}$ --- множество нулей функции $\psi$.

   Обозначим
\begin{equation}\label{prod R}
\Pi(R,z){\hskip -2pt}:\, =\prod\limits_{|z_n| \ge R}E_p\left(\frac{z}{z_n}\right).
\end{equation}
Наша цель на этом эпапе --- оценить  функцию $\Pi(R,z) -1$ в круге
\newline $B(0, a R^{1-\delta})$.


\begin{lemma}\label{lemma_Products}
Пусть $\delta \in (0,1)$, $a>0$ а число $p\in \mathbb N$  подчинено условиям
 $p\geqslant [\rho], \ \alpha{\hskip -4pt}:\, = \delta(p+1) -\rho \geqslant \mu.$
 Тогда  в круге $B\left(0, a R^{1-\delta}\right)$ справедлива оценка
\begin{equation}\label{Products}
 |\Pi(R,z)-1|\le 2C_2 a^{p+1} R^{-\mu}, \qquad C_2 = \frac{2\sigma(p+1)(2e)^\rho}{p+1-\rho},
 \end{equation}
  если только
 \begin{equation}\label{r 2}
 R\geqslant r_2:= \max \left(r_1, \ \left(\frac{a(p+1)}p\right)^{1/\delta},
 \ \left(\frac{C_2 a^{p+1}}{\ln 2}\right)^{1/\mu}\, \right),
 \end{equation}
 где число $r_1$ определено в \eqref{r 1}.
 \end{lemma}

\textbf{Доказательство. } Отметим, что в ходе дальнейшего доказательства эта лемма
будет использоваться при $a=p+1$.
 Имеем \begin{equation*}
\Pi(R,z) - 1 = \exp\left( \sum\limits_{|z_n| \ge R}\ln E_p\left(\frac{z}{z_n}\right)\right) - 1.
\end{equation*}
 Поскольку
 \begin{equation}\label{e w}
 | e^w - 1|\leqslant  |w| e^{|w|}, 
 \end{equation}
 то достаточно оценить сумму
 \begin{equation*}
 w(z)= \sum\limits_{|z_n| \geqslant R}\ln E_p\left(\frac{z}{z_n}\right).
 \end{equation*}
Известно \cite[\S 4.3]{Le}, что при $|\xi|\leqslant
 \frac{p}{p+1}$ справедлива оценка
 \begin{equation*}
 \left|\ln E_p(\xi)\right|\leqslant |\xi|^{p+1}.
 \end{equation*}
 Неравенство  $|z|(p+1) \leqslant p R$ выполнено для всех $z\in \left(0, a R^{1-\delta}\right)$,
 если только $a(p+1) R^{1-\delta} \leqslant p R$ или
 \begin{equation}\label{r_20}
 R\geqslant \left(\frac{a(p+1)}p\right)^{1/\delta}.
 \end{equation}
Следовательно, при таком условии для всех $z\in B\left(0, a R^{1-\delta}\right)$ имеем
 \begin{equation}\label{w}
 |w(z)|\leqslant \sum\limits_{|z_n| \geqslant R}\left| \ln E_p\left(\frac{z}{z_n}\right)\right|\leqslant
  |z|^{p + 1}\sum\limits_{|z_n| \ge R}\frac{1}{|z_n|^{p+1}}.
 \end{equation}
Согласно \eqref{n(r)},  получаем
\begin{multline}\label{prod}
\sum\limits_{|z_n| \ge R}\frac{1}{|z_n|^{p+1}} \leqslant \int_R^\infty\frac{dn(r)}{r^{p+1}}
\leqslant (p+1)\int\limits_R^{\infty}\frac{n(r)}{r^{p+2}}dr\leqslant2(p+1)\sigma(2e)^\rho
\int\limits_R^{\infty}\frac{dr}{r^{p+2-\rho}}\\
\leqslant C_2 R^{\rho- p- 1},
\end{multline}
где  постоянная   $C_2$  определена в \eqref{Products} и $R\geqslant r_1$.
Теперь, полагая $\alpha = \delta(p+1) -\rho$ и учитывая условие $\alpha \geqslant\mu$,
 из оценок \eqref{w} и \eqref{prod} при $z\in B\left(0, a R^{1-\delta}\right)$ получаем
\begin{equation}\label{w 2}
|w(z)| \leqslant |z|^{(p+1)}C_2R^{\rho -p-1} \leqslant C_2 a^{p+1}R^{-\alpha} \leqslant
  C_2 a^{p+1}R^{-\mu}, 
\end{equation}
 если $R>r_1$ и  если  выполнено условие \eqref{r_20}. При дополнительном условии
$C_2 a^{p+1}R^{-\mu} \leqslant \ln 2$ из оценок \eqref{w 2} и \eqref{e w} получаем оценку \eqref{Products}.
Учет всех условий на $R$ приводит к условию \eqref{r 2}. Лемма доказана.

\medskip

 {\it Шаг\ 4.} На  этом этапе мы докажем следующее утверждение.

\begin{lemma}\label{lemma_Polynom}
Обозначим через $W_{kj}$  алгебраические дополнения  элементов $c_{kj} = k^{j-1}$,
находящхся на пересечении  $k$-ой строки и  $j$-го столбца определителя Вандермонда
\begin{equation}\label{W}
W=  \begin{vmatrix}
	1& 1&1&\cdot&\cdot&\cdot& 1\\
 1&2&2^2&\cdot&\cdot&\cdot&2^p\\
 \cdot&\cdot\cdot&\cdot&\cdot&\cdot&\cdot\\
 1&(p+1)&(p+1)^2&\cdot&\cdot&\cdot&(p+1)^p
	\end{vmatrix}	=1!\, 2!\, 3!\ \cdots\, p!.
\end{equation}	
Положим
\begin{equation}\label{A p}
A_p =\frac 1{W} \sum_{k=1}^{p+1}\sum_{j=1}^{p+1} |W_{jk}| k^{-\mu}.
\end{equation}
Пусть $g(z) = a_0+a_1z+ \dots + a_pz^p$ --- многочлен степени $p\geqslant 1$,  причем для всех
$z$  из отрезка $[r, (p+1)r]$  выполнена оценка
$$
\left|e^{g(z)} -1\right| \leqslant C_1|z|^{-\mu}\leqslant \varepsilon, \qquad z \in [r, (p+1)r], \quad\varepsilon{\hskip -3pt}:\, = C_1 r^{-\mu}
 \leqslant \frac 12.
$$
Тогда для всех $z$  в круге $B(0,r)$  выполняется оценка
\begin{equation}\label{Pol-est1}
\left|e^{g(z)} -1\right| \leqslant 2 \varepsilon A_p,
\quad\text{если} \ \, \varepsilon A_p \leqslant 1/4.
\end{equation}
\end{lemma}

{\bf Доказательство}. Положим $h(z)= e^{g(z)} -1$,  тогда
$$
g(z) = \ln(1+ h(z)).
$$
Коэффициент $a_0$  полинома $g$ определяется с точностью до постоянной вида $2\pi k i, \ k\in \mathbb Z.$
Выберем главную ветвь логарифма при $z\in [r, (p+1)r]$, так, чтобы
$$
g(z) = \sum_{k=1}^\infty (-1)^{k-1} \frac {h^k(z)}k.
$$
При $z\in [r, (p+1)r]$  имеем $|h(z)|\leqslant\varepsilon \leqslant 1/2$, поэтому
$$
|g(z)|\leqslant |h(z)| +\frac{|h^2(z)|}2 +\frac{|h^3(z)|}{3} +
\dots \ \, 
\leqslant  |h(z)|(1+|h(z)|).
$$
Следовательно,
\begin{equation}\label{gest}
g(kr) =a_0 +a_1kr + \dots +a_p(kr)^p = \zeta_k, \quad |\zeta_k| \leqslant\varepsilon(1+\varepsilon) k^{-\mu}, \ k=1,2,\dots, p+1.
\end{equation}
Получили систему $p+1$ уравнений, из которой однозначно определяются  коэффициенты $a_j$.
Очевидно, определитель этой системы $\mathcal W$  и алгебраические дополнения $\mathcal W_{jk}$ его элементов,
находящихся на  пересечении $j$-го столбца и $k$-ой строки, равны
$$
\mathcal W = r^{p(p+1)/2} W, \quad  \mathcal W_{kj} = r^{p(p+1)/2 -j+1} W_{kj}.
$$
Согласно правилу Крамера получаем
$$
|a_{j-1}| \leqslant \frac 1W \sum_{k=1}^{p+1}|\zeta_k W_{kj}|\leqslant \varepsilon(1+\varepsilon) r^{-j+1}
\frac 1W \sum_{k=1}^{p+1}|W_{kj}|k^{-\mu}.
$$
Но тогда для всех $z\in B(0, r)$ имеем
$$
|g(z)| \leqslant \frac{\varepsilon(1+\varepsilon)}{W} \sum_{k=1}^{p+1}\sum_{j=1}^{p+1}|W_{kj}|k^{-\mu} =
\varepsilon(1+\varepsilon) A_p.
$$
Теперь воспользуемся оценкой
$$
|e^g -1| \leqslant |g|e^{|g|}\leqslant |g|(1+2|g|), \qquad\text{если} \ \, |g|\leqslant 1/2.
$$
Далее предположим, что
 $\varepsilon A_p \leqslant 1/4$. Ниже в замечании мы покажем, что $A_p > 2^{p+1} -1 \geqslant 3.$
Поэтому такое предположение влечет $\varepsilon <1/12$. Но тогда  для всех $z\in B(0, r)$ получаем
$$
 |e^{g(z)} -1| \leqslant \varepsilon(1+\varepsilon) A_p (1+2\varepsilon(1+\varepsilon) A_p)
 \leqslant 2 \varepsilon A_p.
 $$
Лемма доказана.

\medskip

{\it Шаг 5}. Пусть  обе функции $\psi_1, \psi_2$  принадлежат классу
$\mathcal M$.  Согласно теореме Адамара для них
  справедливы представления
$$
\psi_j(z) = z^{n_j} e^{g_j(z)}\prod_{k=1}^\infty E_p\left(\frac z{z_k^{(j)}}\right), \qquad j=1,2.
$$
где $z_k^{(j)}$ ---  нули  функций $\psi_j$,  а $g_j$ --- многочлены степени $\leqslant p$.
Если нули
функций $\psi_1$  и $\psi_2$  в круге $B(0,R)$ совпадают, то
\begin{equation}\label{represent}
\frac{\psi_2(z)}{\psi_1(z)} = e^{g_2(z) -g_1(z)}\ \frac{\Pi_2(z, R)}{\Pi_1(z, R)},
\end{equation}
где функции $\Pi_1, \Pi_2$ определены \eqref{prod R} и построены по нулям  функций $\psi_1, \psi_2$
соответственно.  Следовательно,
\begin{equation}\label{A}
e^{g_2(z) -g_1(z)} -1 = \left(\frac{\psi_2(z)}{\psi_1(z)} -1 \right) \frac{\Pi_1(z, R)}{\Pi_2(z, R)}
+\left(\frac{\Pi_1(z, R)}{\Pi_2(z, R)}-1\right).
\end{equation}
Оценим правую часть этого равенства  в круге $|z|\leqslant  R^{1-\delta}$.
Далее будут использоваться очевидные неравенства
$$
(1-\eta)^{-1} \leqslant 1+ 3\eta/2, \quad (1+\eta)(1-\eta)^{-1} \leqslant 1+ 3\eta, \ \ \text{если}\ \,  0\leqslant
 \eta \leqslant 1/3.
$$
По условию на положительном луче  справедлива оценка \eqref{2},  поэтому  при $z= r\ge r_0$ получаем
\begin{multline}\label{B}
 \left|\frac{\psi_2(z)}{\psi_1(z)} -1 \right| \leqslant \frac{|\psi_2(z)-1| +|\psi_1(z) -1|}{|\psi_1(z)|}\\
\leqslant \frac{2C_1}{r^\mu}\left( 1-\frac{C_1}{r^\mu}\right)^{-1} \leqslant (2+3\eta_1)\eta_1, \quad\text{если}\ \,
\eta_1{\hskip -1pt}:= \frac{C_1}{r^\mu} \leqslant \frac 13.
\end{multline}
Обозначим $C_3 = 2C_2 (p+1)^{p+1}$.
Согласно лемме 2  при $|z|\leqslant  (p+1) R^{1-\delta}, R\geqslant r_2,$ имеем
\begin{equation}\label{Pi}
\left|\frac{\Pi_1(z, R)}{\Pi_2(z, R)}\right| \leqslant \frac {1+C_3 R^{-\mu}}{1-C_3 R^{-\mu}} \leqslant 1+3\eta_2,
\quad \text{если}\ \, \eta_2 = \frac{C_3}{R^\mu} = \frac{2C_2  (p+1)^{p+1}}{R^\mu}\leqslant \frac 13,
\end{equation}
то есть, если
\begin{equation}\label{r 3}
 R\geqslant r_3{\hskip -1pt}: = \max \left( r_2,  \left(6C_2 (p+1)^{p+1}\right)^{1/\mu}\right).
\end{equation}
Так же, как в \eqref{B},   при $|z|\leqslant (p+1) R^{1-\delta}$  получаем
\begin{equation}\label{C}
\left|\frac{\Pi_1(z, R)}{\Pi_2(z, R)}-1\right| \leqslant \frac {2C_3R^{-\mu}}{1-C_3R^{-\mu}} \leqslant 3\eta_2,
\quad \text{если}\ \, \eta_2{\hskip -1pt}: = \frac{C_3}{R^\mu} \leqslant \frac 13.
\end{equation}

Теперь мы намерены воспользоваться леммой 3. При $ z\in \Delta{\hskip -2pt}: = \left[R^{1-\delta},
(p+1) R^{1-\delta}\right]$ правая часть \eqref{B} принимает
максимальное значение при $z= r= R^{1-\delta}$. Следовательно, для всех $z\in \Delta$  левая часть \eqref{B}
допускает оценку $\leqslant (2+3\eta)\eta$, где
\begin{equation}\label{eta}
\eta{\hskip -2pt}: = \max_{z\in \Delta} \eta_1  =
\frac{C_1}{R^{\mu(1-\delta)}} \leqslant\frac 13.
\end{equation}
Отрезок $\Delta$ входит в круг $ B\left(0, (p+1) R^{1-\delta}\right)$,  поэтому из представления \eqref{A} с учетом оценкок
 \eqref{B}, \eqref{Pi}, \eqref{C},
    при всех $z\in\Delta$ получаем
  \begin{equation*}
\left|e^{g_2(z) -g_1(z)} -1\right|  \leqslant (2+3\eta)\eta(1+3\eta_2)+ 3\eta_2 \leqslant 9\eta,
 \quad \, \text{если} \ \eta_2 \leqslant \eta.
\end{equation*}
Воспользовавшись леммой 3, получим, что для всех
 $z\in B\left(0, а R^{1-\delta}\right)$ справедлива  оценка
\begin{equation}\label{D}
\left|e^{g_2(z) -g_1(z)} -1\right|  \leqslant  18 A_p \eta.
\end{equation}
если $9\eta \leqslant 1/3$  и $9A_p\eta \leqslant 1/4$.
Первое условие  автоматически выполняется, если выполняется второе, так как $A_p > 2^{p+1} - 1\geqslant 3$
(см. ниже замечание 4). Но условие $9A_p\eta \leqslant 1/4$ эквивалентно условию
\begin{equation}\label{r 4}
R\geqslant r_4{\hskip -2pt}: = \left(36\, C_1 A_p \right)^{1/\mu(1-\delta)}.
 \end{equation}
Вспомним еще, что мы использовали неравенство $\eta_2\leqslant \eta$,  которое эквивалентно условию
\begin{equation}\label{r 5}
R\geqslant  r_5{\hskip -2pt}: = \left(2 (p+1)^{p+1} \frac{C_2}{C_1} \right)^{1/\mu\delta}.
\end{equation}
Итак, в круге $B\left(0, (p+1) R^{1-\delta}\right)$ нами получена ключевая оценка
\eqref{D} при условии, что $R\geqslant \max(r_3, r_4, r_5).$

\medskip

{\it Шаг 6.}
Из представления \eqref{represent} получаем
\begin{equation}\label{represent2}
\left|\frac{\psi_2(z)}{\psi_1(z)}-1\right| = \left(1+\left(e^{g_2(z) -g_1(z)}-1\right)\right)\,\left( \frac{\Pi_2(z, R)}{\Pi_1(z, R)}-1\right) + \left(e^{g_2(z) -g_1(z)} -1\right).
\end{equation}
Функция $\Pi_1/\Pi_2 -1$
 оценивается точно также, как  функция  $\Pi_2/\Pi_1 -1$ (изменение порядка индексов значения не имеет).
Поэтому из представления \eqref{represent2}, оценки \eqref{D}  и определения \eqref{eta}
для всех $z\in B\left(0, (p+1) R^{1-\delta}\right)$  получаем
\begin{equation}\label{represent3}
\left|\frac{\psi_2(z)}{\psi_1(z)}-1\right|  \leqslant (1+ 18A_p\eta) 3\eta  + 18 A_p \eta
\leqslant 20\, A_p\, \eta = \frac{ 20\, A_p\,  C_1}{R^{\mu(1-\delta)}}.
\end{equation}
Здесь при переходе ко второму неравенству мы учли, что $18A_p \eta \leqslant 1/2$ и $3\eta \leqslant A_p\eta$.
Этим завершается доказательство теоремы.

\medskip

{\bf 3. Замечания.}
Отметим еще раз, что основой  для настоящей заметки послужила работа
\cite{MaSW}, где для функций вида \eqref{0}  с финитным ядром (то есть, для функций
$K$,  таких, что $K(t)\equiv 0$ при $t\geqslant t_0$) получена оценка
$$
\left|\psi_2(z)/\psi_1(z) - 1\right|\leqslant CR^{-1/3}  \quad \text{при} \ \, |z|\leqslant R^{1/3}, \ \  R\geqslant R_0,
$$
 при условии, что нули функций  $\psi_1, \psi_2$  совпадают в круге радиуса $R$.
 Этот результат получается из теоремы 1 при $p=1$ и $\delta = 2/3$. Дополнительно
 мы имеем в явном виде константы, при которых выполняется это неравенство.

 \medskip

{\bf Замечание 4.}  Константу $A_p$ можно оценить снизу и сверху (в зависимости от $p$).
Мы покажем только, что  $A_p \ > 2^{p+1}-1\geqslant 3$ (эту оценку мы использовали).
Сначала вычислим  определитель $W_{p+1,1}$.  После вынесения из $k$-й строки  этого
определителя множителя $k$ снова получим определитель Вандермонда вида \eqref{W},  но
только размера $p\times p$. Поэтому
$$
W_{p+1,1}\  W^{-1} = 1\cdot 2\cdot 3\cdots p \cdot (p-1)! (p-2)! \cdots 1!\  W^{-1} =1.
$$
Похожим образом вычисляются определители
$$
W_{p-s,1} = (p+1)!(p-s)^{-1} W'_{p-s,1}, \quad s= 0, 1, \dots, p-1,
$$
где  $W'_{p-s,1}$ --- определитель Вандермонда чисел $1,\dots, p-s-1, p-s+1, \dots, p+1$. Он
равен
$$
W'_{p-s,1} = \frac{p!}{s+1}\cdot \frac{(p-1)!}s \cdots \frac {(p-s)!}1\cdot (p-s-2)!\cdots 1!.
$$
После замены $p-s =k$ получаем
$$
\frac{W{_{1k}}}W = \frac {(p+1)!}{k!(p+1-k)!} = C^k_{p+1},
$$
где $ C^k_{p+1}$ --- биномиальные коэффициенты. Но тогда
$$
A_p > \sum_{k=1}^{p+1}C^k_{p+1} = \sum_{k=0}^{p+1}C^k_{p+1} -1 = 2^{p+1} -1.
$$

\medskip

{\bf Замечание 5.}
Функции вида \eqref{0}  ограничены на вещественной оси некоторой постоянной $C$,
зависящей от ядра  $K$. Поэтому совпадение нулей таких функций  в круге $B(0,R)$
влечет оценку их разности на большом отрезке вещественой оси
\begin{multline*}
|\psi_ 2(z) -\psi_1(z)| = |\psi_ 2(z)/\psi_1(z) -1|\ |\psi_1(z)| \leqslant \varepsilon C R^{-(1-\delta)},
\\ \text{если}\ \, z\in\left[- R^{1-\delta}, R^{1-\delta}\right], \ \, R\geqslant R_0.
\end{multline*}

\medskip

{\bf Замечание 6.} Известно \cite[Гл. 2]{Ma}, что производная $K'(t)$  в представлении функции
Йоста \eqref{0} оператора Шредингера экспоненциально убывает, если экспоненциальное убывание имеет
производная потенциала $q$ оператора (порядок убывания сохраняется).  В этом случае
в предствлении \eqref{0}  можно проинтегрировать по частям и получить
$$
\psi(z) = 1 - \frac{K(0)}{iz} - \frac 1{iz} \int_0^\infty K'(t) e^{izt}\, dt.
$$
Если имеются две функции $\psi_1, \psi_2$ такого вида с ядрами $K_1, K_2$,  и известно, что
$K_1(0) = K_2(0)$,  то функции $\tilde\psi_j(z) = \psi_j(z) +K(0)(iz)^{-1}$  в верхней
полуплоскости допускают оценку
$$
|\tilde\psi_j(z)| \leqslant C_1\left( |z| y\right)^{-1}, \quad y=\text{Im}\, z>0.
$$
Поэтому теорема 1 для таких функций может быть применена с показателем $\mu =2$.

Адрес:

МГУ им. М.В.Ломоносова, г.~Москва

email:  valgeynts@gmail.com;  shkalikov@mi.ras.ru


\bigskip

\begin{thebibliography}{99}

\bibitem{Le} B.Ya.Levin. Lectures on entire functions. English revised edition.  Amer. Math. Soc, Providence, RI, 1996

\bibitem{Ma} В.А. Марченко.   Операторы Штурма-Лиувилля и их приложения. Киев: Наукова Думка, 1977.

\bibitem{MaSW} M. Marletta, R. Shterenberg and R. Weikard, On the inverse resonance problem for Schroedinger operators, Comm. Math. Phys., vol. 295, no. 2, 2010, pp. 465-484

\end{thebibliography}
\end{document}